\documentclass[12pt]{article}%
\usepackage{graphicx}
\usepackage[intlimits]{amsmath}
\usepackage{latexsym}
\usepackage{amsfonts}
\usepackage{amssymb}%
\makeatletter
\newfont{\footsc}{cmcsc10 at 8truept}
\newfont{\footbf}{cmbx10 at 8truept}
\newfont{\footrm}{cmr10 at 10truept}
\newtheorem{theorem}{Theorem}

\newenvironment{proof}[1][Proof]{\noindent{\textbf {#1}  }}  {\hfill$\Box$\bigskip}

\begin{document}

\title{Eigenvalues and extremal degrees in graphs}
\author{Vladimir Nikiforov\\Department of Mathematical Sciences, University of Memphis, \\Memphis TN 38152, USA, email: \textit{vnkifrv@memphis.edu}}
\maketitle

\begin{abstract}
Let $G$ be a graph with $n$ vertices, $\mu_{1}\left(  G\right)  \geq...\geq
\mu_{n}\left(  G\right)  $ be the eigenvalues of its adjacency matrix, and
$0=\lambda_{1}\left(  G\right)  \leq...\leq\lambda_{n}\left(  G\right)  $ be
the eigenvalues of its Laplacian$.$ We show that
\[
\delta\left(  G\right)  \leq\mu_{k}\left(  G\right)  +\lambda_{k}\left(
G\right)  \leq\Delta\left(  G\right)  \text{ \ \ for all }1\leq k\leq n,
\]
and%
\[
\mu_{k}\left(  G\right)  +\mu_{n-k+2}\left(  \overline{G}\right)  \geq
\delta\left(  G\right)  -\Delta\left(  G\right)  -1\text{ \ \ for all }2\leq
k\leq n.
\]

Let $\mathcal{G}$ be an infinite family of graphs. We prove that $\mathcal{G}$
is quasi-random if and only if $\mu_{n}\left(  G\right)  +\mu_{n}\left(
\overline{G}\right)  =o\left(  n\right)  $ for every $G\in\mathcal{G}$ of
order $n.$ This also implies that if $\lambda_{n}\left(  G\right)
+\lambda_{n}\left(  \overline{G}\right)  =n+o\left(  n\right)  $ for every
$G\in\mathcal{G}$ of order $n,$ then $\mathcal{G}$ is quasi-random.



\end{abstract}

\section{Introduction}

Our notation is standard (e.g., see \cite{Bol98}, \cite{CDS80}, and
\cite{HoJo88}); in particular, all graphs are defined on the vertex set
$\left\{  1,2,...,n\right\}  ,$ $G\left(  n\right)  $ stands for a graph of
order $n,$ and $\overline{G}$ denotes the complement of $G.$ Writing $A\left(
G\right)  $ for the adjacency matrix of $G$ and $D\left(  G\right)  $ for the
diagonal matrix of its degree sequence, the Laplacian of $G$ is defined as
$L\left(  G\right)  =D\left(  G\right)  -A\left(  G\right)  .$ If $G=G\left(
n\right)  ,$ we order the eigenvalues of $A\left(  G\right)  $ as $\mu
_{1}\left(  G\right)  \geq...\geq\mu_{n}\left(  G\right)  $ and the
eigenvalues of $L\left(  G\right)  $ as $0=\lambda_{1}\left(  G\right)
\leq...\leq\lambda_{n}\left(  G\right)  .$

In this note we prove that if $G=G\left(  n\right)  $ is a graph with minimum
degree $\delta\left(  G\right)  $ and maximum degree $\Delta\left(  G\right)
,$ then%
\begin{equation}
\delta\left(  G\right)  \leq\mu_{k}\left(  G\right)  +\lambda_{k}\left(
G\right)  \leq\Delta\left(  G\right)  \text{ \ \ for all }1\leq k\leq n.
\label{in1}%
\end{equation}
This, in turn, implies that%
\begin{equation}
\mu_{k}\left(  G\right)  +\mu_{n-k+2}\left(  \overline{G}\right)  \geq
\delta\left(  G\right)  -\Delta\left(  G\right)  -1\text{ \ \ for all }2\leq
k\leq n, \label{in2}%
\end{equation}
complementing the well-known inequality $\mu_{k}\left(  G\right)  +\mu
_{n-k+2}\left(  \overline{G}\right)  \leq-1.$

In the second part of this note we give new spectral conditions for
quasi-randomness of graphs. Throughout this note we denote by $\mathcal{G}$ an
infinite family of graphs. Following Chung, Graham, and Wilson \cite{CGW89},
we call a family $\mathcal{G}$ \emph{quasi-random},\ if for every
$G\in\mathcal{G}$ of order $n,$%
\[
\mu_{1}\left(  G\right)  =2e\left(  G\right)  /n+o\left(  n\right)  ,\text{
}\mu_{2}\left(  G\right)  =o\left(  n\right)  ,\text{ and }\mu_{n}\left(
G\right)  =o\left(  n\right)  .
\]
Applying results of \cite{Nik06}, we first prove the following theorem.

\begin{theorem}
\label{th1}A family $\mathcal{G}$ is quasi-random if and only if
\begin{equation}
\mu_{n}\left(  G\right)  +\mu_{n}\left(  \overline{G}\right)  =o\left(
n\right)  \label{con1}%
\end{equation}
for every graph $G\in\mathcal{G}$ of order $n.$
\end{theorem}

This, in turn, implies the following sufficient conditions for
quasi-randomness in terms of Laplacian eigenvalues.

\begin{theorem}
\label{th2}If $\mathcal{G}$ is a family such that
\begin{equation}
\lambda_{n}\left(  G\right)  +\lambda_{n}\left(  \overline{G}\right)
=n+o\left(  n\right)  \label{con2}%
\end{equation}
for every $G\in\mathcal{G}$ of order $n,$ then $\mathcal{G}$ is quasi-random.
\end{theorem}

Since $\lambda_{2}\left(  G\right)  +\lambda_{n}\left(  \overline{G}\right)
=n$ for every $G=G\left(  n\right)  ,$ we also obtain the following theorem.

\begin{theorem}
\label{th3}If $\mathcal{G}$ is a family such that
\[
\lambda_{2}\left(  G\right)  +\lambda_{2}\left(  \overline{G}\right)
=o\left(  n\right)
\]
for every $G\in\mathcal{G}$ of order $n,$ then $\mathcal{G}$ is quasi-random.
\end{theorem}

We leave the extension of the above results to normalized Laplacians to the
interested reader.

\section{Proofs}

\begin{proof}
[\textbf{Proof of inequality (\ref{in1})}]Let $\mathbf{u}_{1},...,\mathbf{u}%
_{n}$ be orthogonal unit eigenvectors to $\lambda_{1},...,\lambda_{n}.$ For
every $k=2,...,n,$ the variational characterization of eigenvalues of
Hermitian matrices (\cite{HoJo88}, p. 178-179) implies that%
\begin{align}
\lambda_{k}\left(  G\right)   &  =\min_{\left\Vert \mathbf{x}\right\Vert
=1,\text{ }\mathbf{x}\bot Span\left\{  \mathbf{u}_{1},...,\mathbf{u}%
_{k-1}\right\}  }\left\langle L\mathbf{x},\mathbf{x}\right\rangle
\label{vc1}\\
\mu_{k}\left(  G\right)   &  =\min_{M\subset\mathbb{R}^{n},\text{ }\dim
M=k-1}\left\{  \max_{\left\Vert \mathbf{x}\right\Vert =1,\text{ }%
\mathbf{x}\bot M}\left\langle A\mathbf{x},\mathbf{x}\right\rangle \right\}
\label{vc2}%
\end{align}
Let $\mathbf{y}$ be such that $\left\langle A\mathbf{y},\mathbf{y}%
\right\rangle $ is maximal subject to $\left\Vert \mathbf{y}\right\Vert =1$
and $\mathbf{y}\bot Span\left\{  \mathbf{u}_{1},...,\mathbf{u}_{k-1}\right\}
$. Letting $\mathbf{y}=\left(  y_{1},...,y_{n}\right)  ,$ we find that
\begin{align*}
\lambda_{k}\left(  G\right)   &  \leq\left\langle L\mathbf{y},\mathbf{y}%
\right\rangle =\sum_{u\in V\left(  G\right)  }d\left(  u\right)  y_{u}%
^{2}-\left\langle A\mathbf{y},\mathbf{y}\right\rangle \leq\Delta\left(
G\right)  -\max_{\left\Vert \mathbf{x}\right\Vert =1,\text{ }\mathbf{x}\bot
Span\left\{  \mathbf{u}_{1},...,\mathbf{u}_{k-1}\right\}  }\left\langle
A\mathbf{x},\mathbf{x}\right\rangle \\
&  \leq\Delta\left(  G\right)  -\min_{M\subset\mathbb{R}^{n},\text{ }\dim
M=k-1}\left\{  \max_{\left\Vert \mathbf{x}\right\Vert =1,\text{ }%
\mathbf{x}\bot M}\left\langle A\mathbf{x},\mathbf{x}\right\rangle \right\}
=\Delta\left(  G\right)  -\mu_{k}\left(  G\right)  ,
\end{align*}
proving the second inequality of (\ref{in1}). The first inequality is deduced
likewise using the dual version of (\ref{vc1}) and (\ref{vc2}).
\end{proof}

\begin{proof}
[\textbf{Proof of inequality (\ref{in2})}]It is known that $\lambda_{k}\left(
G\right)  +\lambda_{n-k+2}\left(  \overline{G}\right)  =n$ for all $2\leq
k\leq n.$ This, in view of (\ref{in1}), implies that
\begin{align*}
n+\mu_{k}\left(  G\right)  +\mu_{n-k+2}\left(  \overline{G}\right)   &
=\lambda_{k}\left(  G\right)  +\lambda_{n-k+2}\left(  \overline{G}\right)
+\mu_{k}\left(  G\right)  +\mu_{n-k+2}\left(  \overline{G}\right) \\
&  \geq\delta\left(  G\right)  +\delta\left(  \overline{G}\right)  \geq
\delta\left(  G\right)  +n-1-\Delta\left(  G\right)  ,
\end{align*}
completing the proof of (\ref{in2}).
\end{proof}

\begin{proof}
[\textbf{Proof of Theorem \ref{th1}}]The necessity of condition (\ref{con1})
is a routine fact, so we shall prove only its sufficiency. Let $G=G\left(
n\right)  ,$ $e\left(  G\right)  =m,$ and set $s\left(  G\right)  =\sum_{u\in
V\left(  G\right)  }\left\vert d\left(  u\right)  -2m/n\right\vert .$ The
following results were obtained in \cite{Nik06}
\begin{align}
\frac{s^{2}\left(  G\right)  }{2n^{2}\sqrt{2m}}  &  \leq\mu_{1}\left(
G\right)  -2m/n\leq\sqrt{s\left(  G\right)  },\label{i1}\\
\mu_{k}\left(  G\right)  +\mu_{n-k+2}\left(  \overline{G}\right)   &
\geq-1-2\sqrt{2s\left(  G\right)  }\text{ \ \ for all }2\leq k\leq
n,\label{i2}\\
\mu_{n}\left(  G\right)  +\mu_{n}\left(  \overline{G}\right)   &  \leq
-1-s^{2}\left(  G\right)  /\left(  2n^{3}\right)  . \label{i3}%
\end{align}
Hence, if (\ref{con1}) holds, (\ref{i3}) implies $\mu_{n}\left(  G\right)
=o\left(  n\right)  ,$ $\mu_{n}\left(  \overline{G}\right)  =o\left(
n\right)  ,$ and $s\left(  G\right)  =o\left(  n^{2}\right)  .$ Thus, from
(\ref{i1}) we obtain $\mu_{1}\left(  G\right)  =2m/n+o\left(  n\right)  .$
Since $\mu_{2}\left(  G\right)  +\mu_{n}\left(  \overline{G}\right)  \leq-1,$
inequality (\ref{i2}) implies that $\mu_{2}\left(  G\right)  =o\left(
n\right)  ,$ completing the proof.
\end{proof}

\begin{proof}
[\textbf{Proof of Theorem \ref{th2}}]According to Grone and Merris
\cite{GrMe94}, $\lambda_{k}\left(  G\right)  \geq\Delta\left(  G\right)  .$
Thus, (\ref{con2}) implies
\[
n-1+\Delta\left(  G\right)  -\delta\left(  G\right)  =\Delta\left(  G\right)
+\Delta\left(  \overline{G}\right)  \leq\lambda_{n}\left(  G\right)
+\lambda_{n}\left(  \overline{G}\right)  =n+o\left(  n\right)  .
\]
Hence,
\[
\Delta\left(  G\right)  -\delta\left(  G\right)  =\Delta\left(  \overline
{G}\right)  -\delta\left(  \overline{G}\right)  =o\left(  n\right)
\]
and (\ref{in1}) implies
\begin{align*}
\mu_{n}\left(  G\right)   &  =-\lambda_{n}\left(  G\right)  +\Delta\left(
G\right)  +o\left(  n\right) \\
\mu_{n}\left(  \overline{G}\right)   &  =-\lambda_{n}\left(  \overline
{G}\right)  +\delta\left(  \overline{G}\right)  +o\left(  n\right)  .
\end{align*}
Adding these two inequalities, in view of (\ref{con2}), we obtain $\mu
_{n}\left(  G\right)  +\mu_{n}\left(  \overline{G}\right)  =o\left(  n\right)
$; the assertion follows from Theorem \ref{th1}.
\end{proof}

\textbf{Acknowledgment }The author is indebted to B\'{e}la Bollob\'{a}s for
his kind support.

\end{document}